\documentclass[letterpaper,10pt, conference,twoside]{ieeeconf}  

\usepackage{amsmath,amssymb,color,cite}
\usepackage{graphicx}

\usepackage{latexsym}
\usepackage{algorithm}
\usepackage{algpseudocode}
\usepackage{verbatim}
\usepackage{cite}
\usepackage[english]{babel}
\usepackage[latin1]{inputenc}
\usepackage{enumerate}

\usepackage{amsthm}

\usepackage{tikz}
\usetikzlibrary{tikzmark}
\usepackage{amsfonts}
\usepackage{amssymb}
\usepackage{amsbsy}
\usepackage{bm}

\DeclareMathOperator{\bcard}{bcard}

\DeclareMathOperator*{\argmin}{arg\,min}

\newcommand{\norm}[1]{\ensuremath{\| #1 \|}}

\let\leq\leqslant
\let\geq\geqslant

\newcommand{\calH}{\ensuremath{\mathcal{H}}}

\newcommand{\calM}{\ensuremath{\mathcal{M}}}
\newcommand{\calN}{\ensuremath{\mathcal{N}}}

\newcounter{todocounter}
\usepackage[colorinlistoftodos]{todonotes}

\newtheorem{theorem}{Theorem}[section]

\newtheorem{corollary}[theorem]{Corollary}
\newtheorem{remark}[theorem]{Remark}

\newtheorem{definition}[theorem]{Definition}

\newtheorem{problem}{Problem}

\title{Informativity for centralized design of distributed controllers
  for networked systems}

\author{Jaap Eising \quad Jorge Cort\'{e}s} 

\begin{document}

\maketitle
\renewcommand{\thefootnote}{\fnsymbol{footnote}}

\footnotetext{This work was partially supported by
	ONR Award N00014-18-1-2828. 
	
	J. Eising and J. Cort\'es are with the Department of Mechanical and Aerospace
	Engineering, University of California, San Diego. \texttt{\{jeising,cortes\}@ucsd.edu}}

\begin{abstract}
  Recent work in data-driven control has led to methods that find
  stabilizing controllers directly from measurements of an unknown
  system. However, for multi-agent systems we are often interested in
  finding controllers that take their distributed nature into
  account. For instance, the full state might not be available for
  feedback at every agent. In order to deal with such information, we
  consider the problem of finding a feedback controller with a given
  block structure based on measured data. Moreover, we provide an
  algorithm that, if it converges, leads to a maximally sparse
  controller.
\end{abstract}

\section{Introduction}

In this paper, we consider the problem of finding distributed
controllers on the basis of measurements of an unknown system. Such
data-driven control problems have garnered a lot of attention
recently, both from the viewpoints of control theory and learning. A
particularly recent development is based on the works by Willems et
al. in \cite{JCW-PR-IM-BLMDM:05} and Markovsky and Rapisarda in
\cite{IM-PR:08}. These works have shifted the focus from the
two-step approach of system identification combined with model based
control towards designing controllers directly from the data.

To be precise, we are interested in finding controllers for
multi-agent systems in the situation where the state matrix is
completely unknown. To compensate for this lack of knowledge, we
assume that we have access to measurements of the input and the
corresponding state collected over a finite time window. In this paper, we take
the viewpoint of the informativity framework, introduced in
\cite{HJVW-JE-HLT-MKC:20}. This means that we find a controller for
the measured system by finding a controller that works for the entire
set of systems consistent with the data. In contrast to
\cite{HJVW-JE-HLT-MKC:20}, we do not assume that the measurements are
exact, but assume that the noise on this time window satisfies bounds
of the form considered in the recent paper
\cite{HJVW-MKC-MM:20}. Among the results of \cite{HJVW-MKC-MM:20} are
conditions that are necessary and sufficient for the problem of
finding a stabilizing
controller. These conditions are
given in the form of the feasibility of linear matrix inequalities
(LMI's), and therefore it is straightforward to check whether they
hold.

However, the controllers found by the aforementioned methods are not
necessarily \textit{distributed}. That is, each agent might require
knowledge of the state of each other agent in order to stabilize the
system. As this might be undesirable or even impossible,
we develop results that take into account the networked structure
of the system. For this, we focus on two different types of
problems. First we consider the problem of designing distributed
controllers according to a given communication graph. That is,
controllers such that agent number $i$ only requires state
measurements from specific other agents. In essence, this requires
us to find state feedback matrices with a given block structure. After
this, we move to the problem of finding controllers with maximal
sparsity. Here we assume that the aforementioned communication graph
is also available for design, and want to find a controller that
guarantees the control objective, yet uses as little communication as
possible.

Our contributions are the following:
\begin{enumerate}
\item We formulate necessary and sufficient conditions, in the form of
  linear matrix inequalities in terms of the data, under which the
  measured system admits a quadratically stabilizing
  controller. These differ from previously known results in the fact
  that we assume $B$ is known.
\item Under certain specific assumptions, we show that the existence
  of such a controller with a given block structure can be checked
  using linear matrix inequalities.
\item We state an algorithm consisting of a repeated convex
  programming problem. If this algorithm converges, we show that it
  finds a controller with maximal sparsity.
\end{enumerate}
Proofs are omitted for space reasons and will appear elsewhere.
\vspace{-8pt}

\subsection*{Literature overview}
As mentioned above, data-driven control has
garnered a lot of attention recently.  Given that it is impossible to
give a complete overview of the field, we refer to the survey paper
\cite{ZH-ZW:13} and the references therein. Some additional work that
needs to be highlighted combines data-driven control and
networks. Specifically, the paper \cite{GB-DSB-FP:21} resolves a number
of data-driven problems regarding complex networks. In
\cite{JJ-HJVW-HLT-MKC-HH:20}, the output synchronization problem is resolved for
leader-follower multi-agent systems. Virtual reference feedback tuning
and $\calH_\infty$ are the topics of \cite{TRVS-ML-PMJVDH:20} and
\cite{TRVS-ML-PMJVDH:21} respectively. Lastly, \cite{JE-HLT:21} provides
conditions on noiseless data for specific analysis problems.

Of course, data-driven control is not only relevant in a context of
networked systems. Many results from more general settings can also be
applied to networks. Some such more recent developments regard the
design of different types of controllers. Specifically, we note the
work on data-driven predictive control
\cite{JC-JL-FD:19,JB-JK-MAM-FA:19,AA-JC:21-csl}, optimal control
\cite{JB-FA:19,CDP-PT:20}, optimization-based control
\cite{GB-MV-JC-EDA:21-tac}, and the behavior-based methods of system level synthesis\cite{AX-NM:21} and input-output parametrizations \cite{LF-BG-AM-GFT}.

Apart from data-driven methods, we should also mention the work on model-based design of distributed or decentralized controllers. First, we note the survey papers \cite{LB:08,LB:14} and the book \cite{DDS:91} and the references therein. A particularly useful method for resolving distributed design problems is provided in the work on quadratic invariance \cite{MR-SL:06,AM-NCM-MCR-SY:12,LL-SL:14}. 

More specifically, finding controllers that are as sparse as possible, while still guaranteeing certain design goals, is also a topic of significant interest. For this topic, a good overview can be found in
\cite{MRJ-NKD:16}. Special mention is made of the paper
\cite{BP-MK-PS:13}, which, like this paper employs LMI's and
\cite{FL-MF-MRJ:13} which deals with an efficient method for resolving
these problems. An important ingredient of most methods noted above is
the idea of reweighted $\ell_1$ minimization of Candes et
al. \cite{EJC-MBW-SPB:08} (see also \cite{SB-LV:04}). Specific
applications of sparse controllers can be found within the field of
power networks \cite{FD-MRJ-MC-FB:14,DKM-FD-SH-SHL-SC-BR-JL:17} and security
\cite{ABA-AT-GJP:21}.

\subsection*{Organization}
The paper is organized as follows. We start with a problem formulation
in Section~\ref{sec:prob}. After this, we introduce the formalities
regarding informativity in Section~\ref{sec:inf}. In particular, that
section focuses on the quadratic stabilizability 
problem, and provides conditions for finding a centralized controller
for each. In Section~\ref{sec:comm} we consider the problem of finding
a controller corresponding to a specific communication graph, which we
resolve for two special cases. We develop an algorithm for finding a
controller that is as sparse as possible in
Section~\ref{sec:sparse}. After this, Section~\ref{sec:sim} illustrates the proposed algorithm using a simulation example. Lastly,
we end the paper with conclusions.

\section{Problem formulation}\label{sec:prob}
Suppose we have a heterogeneous networked system given by $r$ agents
of the form:
\begin{equation}
	x_i(t+1) = \sum_{j=1}^r A_{ij}x_j(t) + B_iu_i(t)
	+w_i(t).
\end{equation}
Denote the state and input dimensions of
agent $i$ by $n_i$ and $m_i$. We can represent the entire system by
\begin{equation}\label{eq:system}
	x(t+1) = A_sx(t)+Bu(t)
	+w(t),
\end{equation}
where 
\[
x(t) = \begin{bmatrix} x_1(t) \\ \vdots \\
	x_r(t)\end{bmatrix}\!\!,\hspace{5pt} u(t) = \begin{bmatrix} u_1(t)
	\\ \vdots \\ u_r(t)\end{bmatrix}\!\!,\hspace{5pt} w(t)
= \begin{bmatrix} w_1(t) \\ \vdots \\ w_r(t)\end{bmatrix}\!\!,
\] 
and
\[
A_s= \begin{bmatrix} A_{11} & \cdots & A_{1r}\\ \vdots & \ddots &
	\vdots \\ A_{r1} & \cdots & A_{rr} \end{bmatrix}\!\!,\hspace{5pt}
B= \begin{bmatrix} B_1 & & 0 \\ & \ddots & \\ 0 & &
	B_r \end{bmatrix}
\!\!. 
\] 
We assume that the input matrix $B$ is known, but that $A_s$ is
unknown. In lieu of this, we assume that we have access to data,
consisting of a finite time window of input and state
measurements. Based on these data we are interested in finding
distributed controllers. In order to formalize this notion, we 
introduce some additional notation.

Suppose that we have a state-feedback controller $K$ that guarantees
some control objective for the system $(A_s,B)$. We can partition $K$
in the same fashion as $A_s$ and $B$, and obtain
\[
K =\begin{bmatrix} K_{11} & \cdots & K_{1r} \\ \vdots & & \vdots
	\\K_{r1} & \cdots & K_{rr}\end{bmatrix},\] with $K_{ij} \in
\mathbb{R}^{p_i\times n_j}$. Note that, if we close the loop, we get
$u(t)=Kx(t)$. In other words, for each agent $i$ we have that
\[
u_i(t) = \sum_{j=1}^r K_{ij}x_j(t). 
\]

An essential observation is the following: If $K_{ij}=0$, then agent
$i$ does not require knowledge of the state of agent $j$ in order to
compute the feedback. As such, we can guarantee the absence of such
dependencies by imposing that certain blocks $K_{ij}$ are equal to
zero. A number of interesting problems now arise.

First of all, there is the problem of \textbf{centralized control},
that is
controller that
stabilizes the system based on measured data. For this
problem we make use the \textit{informativity framework} of
\cite{HJVW-JE-HLT-MKC:20}. This means that we make the observation that we
can only guarantee that a controller attains the objective for the
true system, if it does so for \textit{all systems that could have
	generated the data}.

Following the standard centralized problem, we consider a number
of variants. For the problem of \textbf{control with a given
	communication graph} we suppose that the controller is allowed a
given \textit{communication graph}, that is, for each agent $i$, a set
of `neighboring' agents $\calN_i$ are available for feedback. In line
with the previous discussion, this is equivalent to finding a controller $K$ such that
certain blocks $K_{ij}$ are zero.

Alternatively, we might be tasked with controlling the system as
\textit{efficiently as possible} in a number of ways. The problem of
\textbf{data-driven  control with minimal actuation} consists of finding a
controller with the least number of nonzero block-rows. This
means that the controller acts on the minimal number of
agents. Similarly, we can consider \textbf{data-driven control with minimal
	observation}. By finding a controller with the least number of
nonzero block-columns, the controller is required to measure the state of the least
number of agents.

Lastly, we look at minimizing the number of nonzero blocks in
$K$. We refer to this problem as \textbf{data-driven control with maximal
	sparsity}.

\section{Preliminaries on informativity}\label{sec:inf}
Before we return to the question of distributed controller design, we first formulate results regarding the non-distributed case. In
this section, we use the rather general noise model that was
introduced in \cite{HJVW-MKC-MM:20}. Note that the results presented
here differ from those in the latter paper due to the fact that we
assume $B$ is known.

Suppose that we collect data from system \eqref{eq:system} in the form
of state and input trajectories $x(t)$ and $u(t)$. We capture these
measurements in the matrices:
\begin{align*}
	X &:= \begin{bmatrix} x(0) & \cdots & x(T) \end{bmatrix}, \\ U_-
	&:= \begin{bmatrix} u(0) & \cdots & u(T-1) \end{bmatrix},
\end{align*}
and subsequently write
\begin{align*}
	X_+ &:= \begin{bmatrix} x(1) & \cdots &
		x(T) \end{bmatrix}, \\ X_- &:= \begin{bmatrix} x(0) & \cdots &
		x(T-1) \end{bmatrix}.
\end{align*}
We assume that the noise $w$ is unknown, that is, the samples of
$w(0),w(1),\dots,w(T-1)$ are not measured. However, we do assume that
the noise samples collected in the matrix
\[
W_- := \begin{bmatrix}  w(0) ~ w(1) ~ \cdots ~ w(T-1) \end{bmatrix}
\]
satisfy a given \textit{noise model}. Let
\[
\Phi := \begin{bmatrix}
	\Phi_{11} & \Phi_{12} \\
	\Phi_{12}^\top & \Phi_{22}
\end{bmatrix}  
\]
be such that $\Phi_{11} = \Phi_{11}^\top\in\mathbb{R}^{n\times n}$,
$\Phi_{12}\in\mathbb{R}^{n\times T}$ and $\Phi_{22} = \Phi_{22}^\top
\in\mathbb{R}^{T\times T}$ and $\Phi_{22}< 0$. We now assume that the
noise satisfies
\begin{equation} 
	\label{eq:asnoise}
	\begin{bmatrix}
		I \\ W_-^\top 
	\end{bmatrix}^\top 
	\begin{bmatrix}
		\Phi_{11} & \Phi_{12} \\
		\Phi_{12}^\top & \Phi_{22}
	\end{bmatrix}
	\begin{bmatrix}
		I \\ W_-^\top 
	\end{bmatrix} \geq 0.
\end{equation}

\begin{remark}[Special cases of the noise model]
	{\rm 
		Note that this noise model encompasses, among others, energy bounds
		of the form $W_-W_-^\top\leq Q$, where $Q\in\mathbb{R}^{n\times
			n}$. For a further discussion on the special cases of this noise model,
		we refer to \cite{HJVW-MKC-MM:20}. \hfill$\square$
	}
\end{remark}

Clearly, the true state matrix $A_s$ satisfies 
\[
W_- = X_+ -A_sX_--BU_-,
\]
where $W_-$ satisfies \eqref{eq:asnoise}. As such, it is clear that we
can define the set of all state-matrices compatible with the data as
\[
\Sigma = \{A\in\mathbb{R}^n \mid W_- = X_+-AX_--BU_- \textrm{
	satisfies } \eqref{eq:asnoise}\}. 
\] 
Let $N \in \mathbb{R}^{2n\times 2n}$ be given by:
\begin{equation}\label{eq:def N}
	N := \begin{bmatrix} I & X_+-BU_- \\0& -X_-  \end{bmatrix}\Phi\begin{bmatrix} I& X_+-BU_- \\0& -X_- \end{bmatrix}^\top.
\end{equation}
Then it is straightforward to show that $A\in\Sigma$ if and only if
\begin{equation}\label{eq:qmi sigma}
	\begin{bmatrix} I \\ A^\top\end{bmatrix}^\top N\begin{bmatrix} I \\ A^\top\end{bmatrix} \geq 0.
\end{equation}

As noted, we are interested in determining properties on the true
system, based on the measurements $(U_-,X)$ as described above. Note
that we can only conclude that the true system $(A_s,B)$ has a given
property if $(A,B)$ has that property for all $A\in\Sigma$. This
observation leads to the following definition.

\begin{definition}
	Let $B$ be given. We say that the data $(U_-,X)$ are
	\textit{informative for quadratic stabilization} if there exists a
	feedback gain $K$ and a matrix $P>0$ such that for each
	$A\in\Sigma$:
	\begin{equation}\label{eq:lyap} (A+BK)P(A+BK)^\top
		<P. \end{equation}
\end{definition}


Note that informativity for quadratic stabilization not only requires
all systems in $\Sigma$ to admit the same feedback gain $K$. We also
require all systems in $\Sigma$ to admit the same Lyapunov function
$P$. In particular a shared Lyapunov function is given by $V(x) = x^\top P^{-1} x$. 

We can equivalently write \eqref{eq:lyap} in the form of: 
\begin{equation}\label{eq:qmi stab} 
	\begin{bmatrix} I \\
		A^\top\end{bmatrix}^\top \begin{bmatrix} P-BKPK^\top B^\top & -BKP
		\\ -PK^\top B^\top & -P \end{bmatrix} \begin{bmatrix} I \\
		A^\top\end{bmatrix} > 0.
\end{equation}
This means that characterizing informativity for quadratic
stabilization is equivalent to characterizing when the quadratic
matrix inequality \eqref{eq:qmi sigma} implies \eqref{eq:qmi stab}. 



\begin{theorem}[LMI conditions for stabilization]\label{thm:stab}
	The data $(U_-, X)$	are informative for quadratic stabilization if and only if there
	exist matrices $P>0$, $L$ and scalars $\alpha\geq 0,\beta>0$ such
	that
	\begin{equation}\label{eq:lmi stab}
		\left\lbrack\begin{array}{cc|c}
			P-\beta I&\! 0 &\! BL \\ 0 & 0 & P \\ \hline L^\top B^\top & P
			& P \end{array}\right\rbrack -
		\alpha\left\lbrack\begin{array}{c|c} N & 0 \\ \hline 0&
			0 \end{array}\right\rbrack\geq 0, 
	\end{equation}
	where $N$ is as defined in \eqref{eq:def N}, holds. Moreover, in
	this case the gain $K:=LP^{-1}$ stabilizes all systems in $\Sigma$.
\end{theorem}

\section{Control with a given sparsity structure}\label{sec:comm}
Having resolved the centralized control problems, we move our
attention to distributed controllers. For this, we introduce some
notation.

Let $p\in\mathbb{N}^k$ and $q\in\mathbb{N}^\ell$ such that
$m=\sum_{i=1}^k p_i$ and $n=\sum_{j=1}^\ell q_j$. Given
$M\in\mathbb{R}^{m\times n}$ we can partition it according
to the vectors $p$ and $q$ by
\begin{equation}\label{eq:part}
	M =\begin{bmatrix} M_{11} & \cdots &
		M_{1\ell} \\ \vdots & & \vdots \\ M_{k1} & \cdots &
		M_{k\ell}\end{bmatrix},
\end{equation}
with $M_{ij} \in \mathbb{R}^{p_i\times q_j}$.

We call $\sigma\in\{0,1\}^{k\times \ell}$ a \textit{block sparsity
	structure}, and define the space of matrices corresponding to
$\sigma$ by:
\[ \calM_{p,q}^\sigma := \{ M\in \mathbb{R}^{m\times n} \mid M_{ij}=0 \textrm{ if } \sigma_{ij}=0\}.\] 
As such, it is clear that the problem of control with a given sparsity structure is a special case of the following.

\begin{problem}[Control with a given sparsity
	structure]\label{prob:comm struc}
	Given vectors $p\in\mathbb{N}^k$, $q\in\mathbb{N}^\ell$ such that
	$m=\sum_{i=1}^k p_i$ and $n=\sum_{j=1}^\ell q_j$, and block sparsity
	structure $\sigma\in\{0,1\}^{k\times \ell}$. Provide necessary and
	sufficient conditions for the data $(U_-,X)$ to be informative for
	quadratic stabilization with feedback gain $K\in\calM_{p,q}^\sigma$.
\end{problem}


\begin{remark}[Block partitions and network systems]\label{rem:spec
		cases}
	{\rm 
		Recall the problem formulation of Section~\ref{sec:prob}. There we
		decompose $K$ according to the state and input dimensions of the
		specific subsystems. Clearly, this corresponds to the choice of $r=
		k=\ell$, and the partition $p_i=m_i$ and $q_i=n_i$ for each
		$i$. As such, finding a controller with a given
		communication graph is a special case of Problem~\ref{prob:comm
			struc}.  However, it is important to stress that for
		Problem~\ref{prob:comm struc}, this is not required. An interesting
		alternative case we consider is the case where $r=k$ and
		$p_i=m_i$ for each $i$, but where $\ell=1$. In terms of
		networked systems, this corresponds to actuating only the agents $i$
		for which $\sigma_{i1}=1$. Similarly, we can look at $k=1$, which,
		in terms of the set-up of Section~\ref{sec:prob}, would correspond
		to measuring only agent $i$ for $\sigma_{1i}=1$.  \hfill$\square$}
\end{remark}

Recall that in Theorem~\ref{thm:stab} we formulate conditions for
quadratic stabilization in the form of LMI \eqref{eq:lmi stab} in the
variables $P>0$, $L$, $\alpha\geq 0$ and $\beta>0$. If \eqref{eq:lmi
	stab} is feasible, we can find a suitable feedback gain by taking
$K=LP^{-1}$. However, note that the latter is not linear in the
variables. This means that testing feasibility of the subspace
constraint $LP^{-1} \in\calM_{p,q}^\sigma$ together with the LMI
\eqref{eq:lmi stab} is no longer linear. However, certain special
cases can be resolved in an efficient manner.

First of all, it is straightforward to show that $L$ and $K=LP^{-1}$
have exactly the same (non-)zero rows, regardless of $P$. As such, we
have the following result.

\begin{corollary}[Control with given block-rows]\label{cor:rows} 
	Suppose that $\ell=1$. The data $(U_-,X)$ are informative for quadratic stabilization with
	feedback gain $K\in\calM_{p,q}^\sigma$ if and only if there exists
	$P>0$, $L\in\calM_{p,q}^\sigma$, $\alpha\geq 0$ and $\beta>0$ such
	that \eqref{eq:lmi stab}, where $N$ is defined as in \eqref{eq:def
		N}, holds.
\end{corollary}	

Let $\bar{\sigma}:=I_\ell \in\{0,1\}^{\ell\times \ell}$. Note that if
$P\in\calM_{q,q}^{\bar{\sigma}}$, then $P$ is a block diagonal
$n\times n$ matrix. Moreover, if the matrix
$P\in\calM_{q,q}^{\bar{\sigma}}$ is (block) diagonal, then so is
$P^{-1}$. Furthermore, it is straightforward to prove that if this is
the case, then $L\in\calM_{p,q}^\sigma$ if and only if
$K=LP^{-1}\in\calM_{p,q}^\sigma$.

\begin{remark}[Block diagonal $P$ and networks]
	{\rm 
		Consider the case of networked systems, that is, $\ell=r$ and
		$q_i=n_i$. Then, the assumption that $P$ is block diagonal
		corresponds to the case where
		\[
		x^\top P^{-1} x = \sum_{i=1}^{r} x_i^\top P_{ii}^{-1} x_i.
		\]
		That is, the Lyapunov function is decoupled.\hfill$\square$
	}
\end{remark} 

As such, we can resolve Problem~\ref{prob:comm struc} efficiently
under the additional assumption that $P$ is block diagonal.

\begin{corollary}[Control with diagonal Lyapunov function]\label{cor:blckdiag} 
	The data $(U_-,X)$ are	informative for quadratic stabilization with feedback gain
	$K\in\calM_{p,q}^\sigma$ and Lyapunov matrix
	$0<P\in\calM_{q,q}^{\bar{\sigma}}$ if and only if there exists
	$0<P\in\calM_{q,q}^{\bar{\sigma}}$, $L\in\calM_{p,q}^\sigma$,
	$\alpha\geq 0$ and $\beta>0$ such that \eqref{eq:lmi stab}, where
	$N$ is defined as in \eqref{eq:def N}, holds.
\end{corollary}

\section{Sparse control}\label{sec:sparse}
After considering finding controllers with a given block structure, we
now move to the problem of finding controllers that are as sparse as
possible.

Let $p\in\mathbb{N}^k$ and $q\in\mathbb{N}^\ell$ and let $M$ be a
matrix that is partitioned as in \eqref{eq:part}. We define the
\textit{block cardinality} of a matrix $M$, denoted $\bcard_{p,q} (M)$
as the number of non-zero blocks in $M$. Let $\phi:\mathbb{R}\rightarrow\mathbb{R}$ be the function defined by:  
\[
\phi(x) := \begin{cases} 0 & x=0, \\ 1 &x\neq 0.\end{cases}
\]
Note that the number of nonzero elements of
$\sigma\in\{0,1\}^{k\times \ell}$ is equal to
$\sum_{i=1}^k\sum_{j=1}^\ell\sigma_{ij}$. As such, we have the
following equivalent statements
\[ 
\bcard_{p,q} (M)\! = \!\!
\min\limits_{\sigma \textrm{ s.t. } M\in\calM_{p,q}^\sigma}\sum_{i=1}^k\sum_{j=1}^\ell\sigma_{ij}
\!=\! \sum_{i=1}^k\sum_{j=1}^\ell \phi(\norm{M_{ij}}_F), 
\]
where $\norm{\cdot}_F$ denotes the Frobenius norm.

In the case where $p_i=q_j=1$ for all $i$ and $j$, the block
cardinality is equal to the number of non-zero elements in $M$. This
is often referred to as the $\ell_0$-pseudo norm or simply the
\textit{cardinality} of $M$. It should be stressed, however, that the
(block) cardinality is not a norm, nor a convex function.

This leads us to the following general problem.

\begin{problem}[Control with maximal sparsity]\label{prob:sparse}
	Given vectors $p\in\mathbb{N}^k$, $q\in\mathbb{N}^\ell$ such that
	$m=\sum_{i=1}^k p_i$ and $n=\sum_{j=1}^\ell q_j$ and data $(U_-,X)$
	that are informative for quadratic stabilization, resolve the
	following problem:
	\begin{equation}\label{eq:min prob}\begin{aligned}  
			\textrm{minimize }  \quad  &  \bcard_{p,q} (K), \\
			\textrm{subject to} \quad & \exists P>0 \textup{
				s.t.}\hspace{6pt} \eqref{eq:lyap} \hspace{6pt} \forall
			A\in\Sigma.
	\end{aligned}\end{equation} 	
\end{problem}

\begin{remark}[Interpretation in terms of networked systems]
	{\rm 
		It follows immediately from the reasoning in Remark~\ref{rem:spec
			cases} that this can be used for the problems of control
		with minimal actuation/observation and control with maximal
		sparsity. \hfill$\square$}
\end{remark}

Note that in Problem~\ref{prob:sparse} the objective function is not a
convex function of $K$, and the constraint set is linear in $P$ and
$KP$, but not necessarily in $K$. As such, the problem above is not a
convex problem. This means that the problem can not be resolved by
many standard methods.

An approach that can work for networked systems with a relatively low
number of agents is a simple exhaustive search. In cases where we can
efficiently solve Problem~\ref{prob:comm struc}, we can simply test
feasibility for different block sparsity patterns $\sigma$ with
increasing number of nonzero elements. This method is guaranteed to
provide the correct answer, but scales in a combinatorial way with
$k\ell$.

As a first step towards resolving the minimization
problem~\eqref{eq:min prob}, we formulate the following corollary of
Theorem~\ref{thm:stab}.

\begin{corollary}[Equivalent formulation of control with maximal
	sparsity]\label{cor:alt}
	Resolving
	\eqref{eq:min prob} in Problem~\ref{prob:sparse} is equivalent to:
	\begin{equation}\label{eq:min prob alt}\begin{aligned}  
			\textup{minimize }  \quad  & \bcard_{p,q} (LP^{-1}), \\ 
			\textup{subject to} \quad  & P>0, \exists \alpha\geq 0, \beta>0.\textup { s.t. } \eqref{eq:lmi stab}.
		\end{aligned}
	\end{equation} 	
\end{corollary}

In the following we take an approach based on the method of reweighted
$\ell_1$ minimization, as introduced in~\cite{EJC-MBW-SPB:08}. As such, we
propose a strategy consisting of repeating a weighted optimization
problem and updating the weights, as shown in Algorithm~\ref{alg:LP}.

\begin{algorithm}
	\caption{Reweighted optimization}\label{alg:LP}
	\begin{algorithmic}[1] 
		\State \textbf{Inputs:} Vectors $p\in\mathbb{N}^k$, $q\in\mathbb{N}^\ell$
		with $m=\sum_{i=1}^k p_i$ and $n=\sum_{j=1}^\ell q_j$, matrix $N$ as in \eqref{eq:def N}.
		\State \textbf{Outputs:} $\{(L_{\hat{t}},P_{\hat{t}})\}_{\hat{t}=0}^{t}$ for some $t\geq 1$.
		\State \textbf{Initialize:} Set $t=0$ and find $L_0$ and $P_0>0$ for which there
		exist $\alpha\geq0$ and $\beta>0$ such that~\eqref{eq:lmi stab}
		holds
		\While{$(L_{t-1},P_{t-1})\neq(L_{t},P_{t})$} \smallskip
		\For{$i=1,\ldots k$, $j=1,\ldots, \ell$} \smallskip
		\State Update the weights by:
		\If{$(L_tP_t^{-1})_{ij} \neq 0$} \smallskip
		\State Let $w_{ij}(t) := \dfrac{1}{\norm{(L_tP_t^{-1})_{ij}}_F}$
		\Else \smallskip
		\State Let $w_{ij}(t) := \infty$
		\EndIf\smallskip
		\EndFor\smallskip
		\State Set $f_t(L) := \sum\limits_{i=1}^k\sum\limits_{j=1}^\ell w_{ij}(t)\norm{(LP_t^{-1})_{ij}}_F$
		\State Update the estimates by solving:
		\begin{equation}\label{eq:relaxed}
			\begin{aligned}&(L_{t+1},P_{t+1}):=  \argmin\limits_{(L,P)} f_t(L) , \\ 
				&\textup{subject to } P>0, \exists \alpha\geq 0, \beta>0 \textrm { s.t. } \eqref{eq:lmi stab}
			\end{aligned}
		\end{equation}
		\State Update $t \leftarrow t+1$
		\EndWhile
	\end{algorithmic}
\end{algorithm}

\vspace{-12pt}
Note that the objective function of the optimization problem
\eqref{eq:relaxed} is not dependent on $P$, but on $P_t$. As such, it
is straightforward to show that the objective function is a convex
function of $L$. Furthermore, the constraint set is given by an LMI,
making it straightforward to resolve \eqref{eq:relaxed}.

\begin{theorem}[If reweighted optimization converges, its output solves the stabilization problem with maximal sparsity]\label{thm:alg}
	Given vectors $p\in\mathbb{N}^k$, $q\in\mathbb{N}^\ell$ such that
	$m=\sum_{i=1}^k p_i$ and $n=\sum_{j=1}^\ell q_j$. Suppose that the data $(U_-,X)$ are informative
	for quadratic stabilization. Then, we can initialize
	Algorithm~\ref{alg:LP}. Moreover, if $(L_{t-1},P_{t-1}) =
	(L_{t},P_{t})$, and we denote $L:=L_{t}$ and $P:=P_{t}$ then
	$LP^{-1}$ is the minimizer of \eqref{eq:min prob}.
\end{theorem}
\vspace{-6pt}
It is important to realize that Theorem~\ref{thm:alg} only gives
sufficient conditions for resolving Problem~\ref{prob:sparse}, since
we have not formulated conditions under which the algorithm
converges.

\vspace{-10pt}

\section{Simulations}\label{sec:sim}
Let the true system be given by $3$ agents, where $n_i=2$, $m_i=1$ and
$B_i=\begin{bmatrix} 1\\ 0\end{bmatrix}$ for each $i\leq 3$. Assume
that the true state matrix is given by:
\[
A_s= \frac{3}{5} \begin{bmatrix}
	1&0&1&0&0&0\\
	1&1&1&1&0&0\\
	0&0&0&0&1&0\\
	0&0&0&0&1&1\\
	1&0&0&0&1&0\\
	1&1&0&0&1&1\end{bmatrix}, \textrm{ and } B
= \begin{bmatrix}1&0&0\\0&0&0\\0&1&0\\0&0&0\\0&0&1\\0&0&0\end{bmatrix}.
\]
We generate measurements using Matlab, by choosing an initial
condition $x(0)$, inputs $u(t)$ and noise $w(t)$ randomly for
$t=0,\ldots 9$, such that $W_-W_-^\top\leq \tfrac{1}{20}I$. The
precise measurements can be found in \eqref{eq:measurements}. 

We use Yalmip \cite{JL:04} with Mosek as a solver in
combination with Theorem~\ref{thm:stab} to resolve the informativity
problem. The
solver returns $P>0$, $L$, $\alpha\geq 0$ and $\beta>0$ such that LMI
\eqref{eq:lmi stab} holds. As such, the data are informative for
quadratic stabilization. In addition this results in the stabilizing
feedback gain $K_1=LP^{-1}$ for all $A\in\Sigma$, given by:
\[
\scriptsize \begin{bmatrix} 
	-0.071335\!\!\!&0.53919\!\!\!&-0.36814\!\!\!&0.23887\!\!\!&-0.72051\!\!\!&-0.74332\\
	0.088392\!\!\!&0.091179\!\!\!&-0.38196\!\!\!&-0.37764\!\!\!&-0.64738\!\!\!&-0.060889\\
	-0.076069\!\!\!&0.54351\!\!\!&0.11392\!\!\!&0.10647\!\!\!&-1.2478\!\!\!&-0.66924
\end{bmatrix}
\]
It can be easily verified that this gain indeed stabilizes the true
system. However, since $K$ is a full matrix we see that in order to compute the input $u_i$, we require for each $j=1,\ldots,r$ the state $x_j$. That is, the controller is not sparse. As such, we move
our attention to Problem~\ref{prob:sparse}, the problem of control
with maximal sparsity. Note that we are not in the situation
Corollary~\ref{cor:rows} or Corollary~\ref{cor:blckdiag}. As such, we
have no efficient way of resolving Problem~\ref{prob:comm struc}. This
prevents us from performing an exhaustive search for a maximally
sparse controller. In addition, note that finding a feedback gain
with, for example, less than or equal to 4 nonzero blocks would
require us to check up to $255$ different sparsity patterns.
\begin{figure*}[t]
	\begin{equation}\label{eq:measurements}\small\begin{array}{rrl}
			\vspace{3pt}X=& &\left\lbrack\begin{array}{lllllllllll}
				0.75274\!&\!1.2276\!&\!1.5028\!&\!1.4546\!&\!2.2505\!&\!3.2402\!&\!4.0554\!&\!4.8123\!&\!5.0687\!&\!5.8844\!&\!7.3989 \\
				0.48475\!&\!1.6001\!&\!2.0504\!&\!3.067\!&\!5.4602\!&\!8.2031\!&\!11.736\!&\!16.6432\!&\!22.2254\!&\!28.3431\!&\!36.5077 \\
				0.62701\!&\!0.28679\!&\!0.56613\!&\!1.9483\!&\!1.9467\!&\!2.3218\!&\!3.3144\!&\!3.4338\!&\!3.7058\!&\!5.0454\!&\!5.0957 \\
				0.80199\!&\!0.30132\!&\!0.99168\!&\!2.6294\!&\!4.0138\!&\!5.7936\!&\!8.6326\!&\!12.1519\!&\!16.2375\!&\!21.5731\!&\!27.317 \\
				0.11059\!&\!0.60892\!&\!1.6934\!&\!1.9273\!&\!2.9284\!&\!3.9676\!&\!4.7534\!&\!5.4348\!&\!6.8431\!&\!7.5784\!&\!8.6763 \\
				0.39059\!&\!1.0436\!&\!2.6886\!&\!4.7617\!&\!6.7268\!&\!10.4199\!&\!15.4993\!&\!21.6268\!&\!29.1105\!&\!37.9496\!&\!47.8538
			\end{array}\right\rbrack \\
			\vspace{3pt} U_-=& &\left\lbrack \begin{array}{lllllllllll}
				0.39914\!&\!0.59328\!&\!0.21324\!&\!0.20845\!&\!0.72101\!&\!0.71757\!&\!0.39015\!&\!0.12077\!&\!0.61899\!&\!0.8402\\
				0.22042\!&\!0.20061\!&\!0.93207\!&\!0.79012\!&\!0.56395\!&\!0.93289\!&\!0.58158\!&\!0.4449\!&\!0.9393\!&\!0.54806\\
				0.090819\!&\!0.59133\!&\!0.0087293\!&\!0.89861\!&\!0.85981\!&\!0.42837\!&\!0.14863\!&\!0.69451\!&\!0.43057\!&\!0.59851
			\end{array}\right\rbrack \\
			W_-= &\!\!\!\!\!\!10^{-4}\!\!\!\!\!&\left\lbrack \begin{array}{llllllllll} 
				0.56402\!&\!0.85894\!&\!0.078075\!&\!0.30536\!&\!0.81527\!&\!0.68118\!&\!0.19788\!&\!0.30939\!&\!0.66536\!&\!0.8844\\
				0.21199\!&\!0.93952\!&\!0.38109\!&\!0.63732\!&\!0.34066\!&\!0.82892\!&\!0.067992\!&\!0.74664\!&\!0.63701\!&\!0.16617\\
				0.020618\!&\!0.17608\!&\!0.26612\!&\!0.25169\!&\!0.81665\!&\!0.99683\!&\!0.21282\!&\!0.0048493\!&\!0.20266\!&\!0.57528\\
				0.61413\!&\!0.1923\!&\!0.19338\!&\!0.42205\!&\!0.42013\!&\!0.11501\!&\!0.24711\!&\!0.46404\!&\!0.91496\!&\!0.25192\\
				0.10097\!&\!0.13537\!&\!0.88955\!&\!0.63512\!&\!0.39169\!&\!0.35093\!&\!0.91207\!&\!0.34179\!&\!0.69204\!&\!0.14824\\
				0.35514\!&\!0.51728\!&\!0.61431\!&\!0.50191\!&\!0.33043\!&\!0.84755\!&\!0.31911\!&\!0.24342\!&\!0.93888\!&\!0.53028
			\end{array}\right\rbrack\end{array}\end{equation}
	\hrulefill
	\vspace{-20pt}
\end{figure*}
In line with Section~\ref{sec:sparse}, we implement
Algorithm~\ref{alg:LP} numerically. This requires making a number of
straightforward changes regarding machine precision to the pseudo
code. Again, we apply Yalmip with the solver Mosek. After 21
iterations, the algorithm has stabilized up to the required
precision. The corresponding feedback gain, denoted $K_{21}$, is found as:
\[\vspace{0.5em}
\scriptsize \begin{bmatrix}
	0\!&\!0\!&\!0\!&\!0\!&\!0\!&\!0\\
	0\!&\!0\!&\!-0.19134\!&\!-0.048629\!&\!0\!&\!0\\
	-0.94584\!&\!-0.052014\!&\!-0.11946\!&\!0.073348\!&\!-0.98268\!&\!-0.14899
\end{bmatrix} \]
As such, we have obtained a feedback gain with just 4 nonzero blocks
that stabilizes all systems in $\Sigma$.

\section{Conclusions}
We have considered data-driven distributed and sparse control. In
particular, we started with defining and resolving informativity
problems regarding centralized stabilization.
As such, we formulated conditions under which a controller guarantees
stabilization 
for all systems compatible with
given measurements. After this, we have considered the same problem
while restricting the allowed controllers to those corresponding to a
given communication graph. For two specific cases, it was shown that
efficient solutions are possible. Lastly we formulated an algorithm
whose steps can be calculated efficiently. If this algorithm
converges, it results in the most sparse stabilizing controller for
all systems compatible with the data.  Future work will investigate
the synthesis of stabilizing controllers with a given sparsity
structure (cf. Problem~\ref{prob:comm struc}) for the general case,
the application of efficient solution methods to the stabilization
with maximal sparsity (cf. Problem~\ref{prob:sparse}), establishing
convergence of Algorithm~\ref{alg:LP}.
A last problem of interest is the case
where a block structure of the state matrix is known, in addition to
one for the controller. Employing such knowledge of the network structure of the system Being able to use this knowledge, might lead to stronger results.

\bibliography{../bib/alias,../bib/JC,../bib/Main,../bib/Main-add,../bib/New,../bib/FB}
\bibliographystyle{IEEEtran}

\end{document}